\documentclass[12pt,oneside]{article}
\usepackage{amsmath,amssymb,amsfonts,amsthm}
\usepackage{stackrel}
\newcommand{\set}[1]{\left\{#1\right\}}

\newcommand\overcirc[1]{\raisebox{10pt}{\tiny$\circ$}{\kern-6.5pt}\mbox{$#1$}}
\newcommand\undersym[2]{\raisebox{-6pt}{\tiny$#2$}{\kern-5pt}\mbox{$#1$}}

\def\Section#1{\vspace{30truept}\addtocounter{section}{1}\setcounter{thm}{0}
\setcounter{equation}{0}{\noindent\Large\bf
    \arabic{section}.~~#1}\par \vspace{12pt}}

\newtheorem{thm}{Theorem}[section]

\newtheorem{lem}[thm]{Lemma}
\newtheorem{prop}[thm]{Proposition}
\newtheorem{defn}[thm]{Definition}

\numberwithin{equation}{section}

\begin{document}
\title{\bf{New Conformal Invariants in Absolute Parallelism Geometry}}
\author{\textbf{Nabil L. Youssef$^{\,1}$, A. Soleiman$^{2}$ and Ebtsam H. Taha$^{1}$}}
\date{}
\maketitle                     
\vspace{-1.16cm}
\begin{center}
{$^{1}$Department of Mathematics, Faculty of Science, Cairo
University, Giza, Egypt}
\end{center}
\vspace{-20pt}
\begin{center}
{$^{2}$Department of Mathematics, Faculty of Science, Benha
University, Benha,
 Egypt}
\end{center}

\begin{center}
nlyoussef@sci.cu.edu.eg,\, nlyoussef2003@yahoo.fr\\ amr.hassan@fsci.bu.edu.eg,\, amrsoleiman@yahoo.com\\ebtsam.taha@sci.cu.edu.eg,\, ebtsam.h.taha@hotmail.com
\end{center}

\vspace{0.1cm}
\maketitle
\begin{flushright}
\emph{ Dedicated to the meomery of Waleed A. Elsayed}\\
\end{flushright}

\noindent{\bf Abstract.} The aim of the present paper is to investigate conformal changes in absolute parallelism geometry. We find out some new conformal invariants in terms of the Weitzenb\"ock connection and the Levi-Civita connection of an absolute parallelism space.




\Section{Introduction}

Conformal transformations play an important role not only in differential geometry but also in application to other branches of science, especially in physics. Conformal changes have been investigated in depth in Riemannian geometry by many authors, cf. for example \cite{Knebelman, Obata, Yan.1}. They have been also investigated  thoroughly in Finsler geometry  \cite{Abed, Hashiguchi, Izumi1, Izumi2, Matsumoto, NAS1, NAS2}.

In the present paper we investigate the conformal changes in absolute parallelism (AP-) geometry.
AP-geometry \cite{Br.1, Rob.1,  Wanas, Wan.2, Waleed, Sid-Ahmed} has a very wide spectrum of applications in physics, especially in the geometrization of physical theories such general relativity and geometric field theories \cite{WN, NW, Shirafuji, charge, WYS}. It is for this reason that we have been motivated to study conformality in such geometry.

We investigate the conformal change of the canonical connections of the space as well as their most important associated tensors. We find out some new conformal invariants related to some conformal connections in AP-geometry. These invariants seem to be promising for applications.


\newpage
\Section{A brief account of AP-geometry}

In this section, we give a brief account of the geometry of parallelizable manifolds or absolute parallelism geometry. For more details, we refer for example to \cite{Br.1, Rob.1, Wan.2, Waleed, Sid-Ahmed}.

A parallelizable manifold \cite{Br.1} is an n-dimensional smooth manifold $M$ which admits $n$
independent global vector fields  $\stackrel[i]{}{\lambda}   (i = 1,...,n)$ on $M$.
Such a space is also known in the literature as an absolute parallelism (AP-) space or a teleparallel space.

  Let $\stackrel[i]{}{\lambda}^{\mu}   (\mu= 1,...,n)$ be the coordinate components of the $i$-th vector field $\stackrel[i]{}{\lambda}$. The covariant components of ${\lambda}^{\mu}$ are given via the relations
\begin{equation}\label{eq.1}
 \stackrel[i]{}{\lambda}^{\mu}\stackrel[i]{}{\lambda}\,_{\!\!\!\nu}=\delta^{\mu}_{\nu}, \quad
 \stackrel[i]{}{\lambda}^{\mu}\stackrel[j]{}{\lambda}\,_{\!\!\!\mu}=\delta_{ij}   .
\end{equation}
 The Einestein summation convention is applied on both Latin (mesh) and Greek (world) indices, where all Latin indices are written downward.

The $n^3$ functions $\Gamma ^{\alpha}_{\mu \nu}$ defined by
\begin{equation}\label{eq.2}
 \Gamma ^{\alpha}_{\mu \nu}:= \stackrel[i]{}{\lambda}^{\alpha} \stackrel[i]{}{\lambda} \,_{\!\!\!\mu,\nu}
\end{equation}
 transform as the coefficients of a linear connection under a change of coordinates, where the comma denotes partial differentiation with respect to the coordinate functions. The connection $\Gamma ^{\alpha}_{\mu \nu}$ is known as the Weitzenb\"ock connection. It is clearly non-symmetric; let its  torsion tensor be denoted by $\Lambda ^{\alpha}_{\mu \nu}:=\Gamma ^{\alpha}_{\mu \nu}-\Gamma ^{\alpha}_{\nu \mu}$. One can show that  $\lambda_{\mu\mid\nu}=0={\lambda^{\mu}}_{\mid\nu},$
  where the stroke denotes covariant differentiation with respect to the Weitzenb\"ock connection (\ref{eq.2}).
  This relation is known in the literature as the AP-condition.
  It is to be noted that the curvature tensor of the Weitzenb\"ock  connection vanishes identically.

  The parallelization vector fields $\stackrel[i]{}{\lambda}$ define a Riemannian metric on $M$ given by
  $$ g_{\mu\nu}:=\stackrel[i]{}{\lambda}\,_{\!\!\!\mu}\stackrel[i]{}{\lambda}\,_{\!\!\!\nu}\nonumber$$
  with inverse
$ g^{\mu\nu}:=\stackrel[i]{}{\lambda}^{\mu}\stackrel[i]{}{\lambda}^{\nu}$. The Levi-Civita connection associated with $g_{\mu\nu}$ is given by the Christoffel symbols:
     \begin{equation*}\label{eq.5}
    \overcirc{\Gamma}^\alpha_{\mu\nu}:=\frac{1}{2}\,g^{\alpha\epsilon}
    \set{g_{\epsilon\nu,\mu}+g_{\epsilon\mu,\nu}-g_{\mu\nu,\epsilon}}.
   \end{equation*}
The contortion tensor $\gamma^\alpha_{\mu\nu}$ is given by
$$\gamma^\alpha_{\mu\nu}:=\Gamma ^{\alpha}_{\mu \nu}- \overcirc{\Gamma}^\alpha_{\mu\nu}= \undersym{\lambda}{i}^{\alpha}\,\,
\undersym{\lambda}{i}_{\mu\,;\,\nu},$$
where the semicolon denotes covariant differentiation with respect to \, $\overcirc{\Gamma^\alpha_{\mu\nu}}$.

Summing up, there are associated to an AP-space two remarkable natural connections: the Weitzenb\"ock connection $\Gamma ^{\alpha}_{\mu \nu}$ and the Levi-Civita connection\, $\overcirc{\Gamma}^{\alpha}_{\mu\nu}$.


\newpage
\Section{Conformal change of AP-space}

In this section, we investigate the conformal change of the natural connections, defined
 in an AP-space, as well as their associated tensors.

Let $(M,\,\undersym{\lambda}{i})$ be an $n$-dimensional AP-space.
Let $\Gamma^{\alpha}_{\mu\nu}$ and $\overcirc{\Gamma}^{\alpha}_{\mu\nu}$ be the Weitzenb\"ock and Levi-Civita connections, respectively.

\begin{defn}\label{def.1}
  Two AP-spaces $(M,\stackrel[i]{}{\lambda})$ and $(M,\overline{\stackrel[i]{}{\lambda}})$  are said to be
conformal (or conformally related) if there exists a positive smooth function
$\rho(x)$  such that
\begin{equation}\label{eq.7}
 \stackrel[i]{}{\overline{\lambda}}^{\,\mu}=e^{-\rho(x)}\stackrel[i]{}{{\lambda}}^{\mu}  \,\,( \text{or} \,
  \stackrel[i]{}{\overline{\lambda}}\,_{\!\!\!\mu}=e^{\rho(x)}\stackrel[i]{}{{\lambda}}\,_{\!\!\!\mu}),
\end{equation}
or, equivalently,
 $$\overline{g}_{\mu\nu}=e^{2\rho(x)}{g}_{\mu\nu}.$$
 \end{defn}

\vspace{5pt}
Now, we present the conformal change of the most important geometric objects associated with an AP-space. Let $\rho_\mu:=\frac{\partial\rho}{\partial x^\mu}$ and $\rho^\mu:=g^{\mu\nu}\rho_\nu$.
\begin{prop}\label{th.1a}
Under the conformal change $(\ref{eq.7})$, we have:
\begin{description}
  \item[(a)] The  Weitzenb\"ock connections $\Gamma^{\alpha}_{\mu\nu}$ and $\overline{\Gamma}^{\,\alpha}_{\mu\nu}$ are related by
  $$ \overline{\Gamma}^{\,\alpha}_{\mu\nu}=\Gamma^{\alpha}_{\mu\nu}+\delta^{\alpha}_{\mu}\rho_{\nu}.$$
  \item[(b)] The  torsion tensors  $\Lambda^{\alpha}_{\mu\nu}$ and $\overline{\Lambda}^{\,\alpha}_{\mu\nu}$ of\, $\Gamma^{\alpha}_{\mu\nu}$ and $\overline{\Gamma}^{\,\alpha}_{\mu\nu}$  are related by
  $$ \overline{\Lambda}^{\,\alpha}_{\mu\nu}=\Lambda^{\alpha}_{\mu\nu}+
  (\delta^{\alpha}_{\mu}\rho_{\nu}-\delta^{\alpha}_{\nu}\rho_{\mu}).$$
 \end{description}
\end{prop}

\begin{prop}\label{th.4}
Under the conformal change $(\ref{eq.7})$, we have:
\begin{description}
  \item[(a)] The Levi-Civita connections\,\, $\overcirc{\Gamma}^{\alpha}_{\mu\nu}$
  and\, $\overline{\overcirc{\Gamma}}\,^{\alpha}_{\mu\nu}$ are related by
  $$ \overline{\overcirc{\Gamma}}\,^{\alpha}_{\mu\nu}
  =\overcirc{\Gamma}^{\alpha}_{\mu\nu}
  +(\delta^{\alpha}_{\mu}\rho_{\nu}+\delta^{\alpha}_{\nu}\rho_{\mu}-g_{\mu\nu}\rho^{\alpha}),$$

    \item[(b)] The  curvature tensors\,  $\overcirc{R}\,^{\alpha}_{\mu\nu\sigma}$ and\,  $\overline{\overcirc{R}}\,^{\alpha}_{\!\mu\nu\sigma}$ of\,\, $\overcirc{\Gamma}\,^{\alpha}_{\mu\nu}$
  and\, $\overline{\overcirc{\Gamma}}\,^{\alpha}_{\mu\nu}$ are related by\footnote{By the symbol $\mathfrak{U}_{\mu\nu}$ we mean:\, $\mathfrak{U}_{\mu\nu}\{A_{\mu\nu}\}:=A_{\mu\nu}-A_{\nu\mu}$}
  $$ \overline{\overcirc{R}}\,^{\alpha}_{\mu\nu\sigma}=\overcirc{R}\,^{\alpha}_{\mu\nu\sigma}
  + \mathfrak{U}_{\nu\sigma}  \set{\delta^{\alpha}_{\sigma}S_{\mu\nu}-g_{\mu \sigma} S^{\alpha}_{\nu}},$$
where $S_{\mu \nu}:= \rho_{\mu; \nu}-\rho_{\mu}\rho_{\nu}-\frac{1}{2}g_{\mu\nu}\rho^{2}$,\,\,
$\rho^{2}:=\rho^{\epsilon}\rho_{\epsilon}$ and\,
 $S^{\alpha}_{\nu}:=g^{\alpha\epsilon} S_{\epsilon\nu}$.\\
 \item[(c)] The contortion tensors $\gamma^\alpha_{\mu\nu}$ and $\overline{\gamma}^\alpha_{\mu\nu}$ are related by
$$\overline{\gamma}^\alpha_{\mu\nu}=\gamma^\alpha_{\mu\nu} -\delta^\alpha_\nu\rho_\mu+g_{\mu\nu}\rho^\alpha.$$
  \end{description}
 \end{prop}


\Section{Conformal connections and invariant tensors}

In this section, we construct three conformally invariant tensors in an AP-space.

\begin{defn} Let $(M,\stackrel[i]{}{\lambda})$ be  an AP-space. A linear connection $\Omega^{\alpha}_{\mu\nu}$ is said
to be conformal if it is conformally invariant under the conformal change $(\ref{eq.7})$:  $\overline{\Omega}^{\,\alpha}_{\mu\nu}=\Omega^{\alpha}_{\mu\nu}$.
\end{defn}

\vspace{5pt}

 \noindent\textbf{Theorem A.}\label{th.conf1}
\emph{Let $(M,\stackrel[i]{}{\lambda})$ be  an AP-space of dimension $n\geq2$. Let $\Lambda^{\alpha}_{\mu\nu}$
and $C_{\mu}:=\Lambda^{\epsilon}_{\epsilon\mu}$ be respectively the torsion and the contracted torsion associated with the Weitzanb\"ock connection $\Gamma^{\alpha}_{\mu\nu}$.
The tensors
\begin{eqnarray*}
  T^{\alpha}_{\mu\nu}&:=&{\Lambda}^{\alpha}_{\mu\nu}-\frac{1}{(n-1)}
\set{\delta^{\alpha}_{\mu}C_{\nu}-\delta^{\alpha}_{\nu}C_{\mu}}, \\
  K^{\alpha}_{\mu\nu\sigma}&:=&\frac{1}{(n-1)}
\set{\delta^{\alpha}_{\mu}C_{\nu,\sigma}-\delta^{\alpha}_{\mu}C_{\sigma,\nu}},
\end{eqnarray*}
are conformally invariant. Moreover, the tensors $T^{\alpha}_{\mu\nu}$ and $K^{\alpha}_{\mu\nu\sigma}$
are precisely the torsion and curvature tensors of a conformal connection on $M$.}

 \begin{proof} Under the conformal change (\ref{eq.7}), using Proposition \ref{th.1a}\textbf{(b)}, we have
 \begin{equation}\label{C}
   \overline{C}_{\nu}=C_{\nu}+(n-1)\rho_{\nu}
 \end{equation}
 From which together with Proposition \ref{th.1a}\textbf{(b)}, we obtain
 \begin{eqnarray*}
   {\overline{T}}^{\:\alpha}_{\mu\nu}&=&{\overline{\Lambda}}^{\:\alpha}_{\mu\nu}-\frac{1}{(n-1)}
\set{\delta^{\alpha}_{\mu}\overline{C}_{\nu}-\delta^{\alpha}_{\nu}\overline{C}_{\mu}} \\
   &=&\Lambda^{\alpha}_{\mu\nu}+  (\delta^{\alpha}_{\mu}\rho_{\nu}-\delta^{\alpha}_{\nu}\rho_{\mu})
   -\frac{1}{(n-1)}\set{\delta^{\alpha}_{\mu}(C_{\nu}+(n-1)\rho_{\nu})-\delta^{\alpha}_{\nu}(C_{\mu}+(n-1)\rho_{\mu})}, \\
&=&{\Lambda}^{\alpha}_{\mu\nu}-\frac{1}{(n-1)}
\set{\delta^{\alpha}_{\mu}C_{\nu}-\delta^{\alpha}_{\nu}C_{\mu}} = T^{\alpha}_{\mu\nu}.
 \end{eqnarray*}
Similarly, from (\ref{C}), we get
 \begin{eqnarray*}
 {\overline{K}}^{\:\alpha}_{\mu\nu\sigma}&=&\frac{1}{(n-1)}
\set{\delta^{\alpha}_{\mu}\overline{C}_{\nu,\sigma}-\delta^{\alpha}_{\mu}\overline{C}_{\sigma,\nu}}\\
&=&\frac{1}{(n-1)}\set{\delta^{\alpha}_{\mu}(C_{\nu}+(n-1)\rho_{\nu})_{\!,\sigma}-\delta^{\alpha}_{\mu}(C_{\sigma}+(n-1)\rho_{\sigma})_{\!,\nu}}\\
&=&\frac{1}{(n-1)}\set{\delta^{\alpha}_{\mu}C_{\nu,\sigma}-\delta^{\alpha}_{\mu}C_{\sigma,\nu}}=
K^{\alpha}_{\mu\nu\sigma}.
 \end{eqnarray*}
This means that  $T^{\alpha}_{\mu\nu}$  and $K^{\:\alpha}_{\mu\nu\sigma}$ are  conformally invariant.

\par Now, let us define the connection:
\begin{equation}\label{e.9}
\mathbf{\Gamma}^{\alpha}_{\mu\nu}:={\Gamma}^{\alpha}_{\mu\nu}-\frac{1}{(n-1)}\delta^{\alpha}_{\mu}C_{\nu}.
\end{equation}
From (\ref{e.9}) and (\ref{C}) together with Proposition \ref{th.1a}\textbf{(a)}, we conclude that
\begin{eqnarray*}
 \mathbf{\overline{\Gamma}}^{\alpha}_{\mu\nu}&=&{\overline{\Gamma}}^{\alpha}_{\mu\nu}-\frac{1}{(n-1)}\delta^{\alpha}_{\mu}\overline{C}_{\nu}\\
 &=& \Gamma^{\alpha}_{\mu\nu}+\delta^{\alpha}_{\mu}\rho_{\nu}-\frac{1}{(n-1)}\delta^{\alpha}_{\mu}(C_{\nu}+(n-1)\rho_{\nu})\\
 &=&{\Gamma}^{\alpha}_{\mu\nu}-\frac{1}{(n-1)}\delta^{\alpha}_{\mu}C_{\nu}=
 \mathbf{\Gamma}^{\alpha}_{\mu\nu},
  \end{eqnarray*}
  which shows that the connection $\mathbf{\Gamma}^{\alpha}_{\mu\nu}$ is conformal.
Next, we prove that the tensors $T^{\alpha}_{\mu\nu}$ and $K^{\alpha}_{\mu\nu\sigma}$
are respectively the torsion and the curvature tensors of the conformal connection $\mathbf{\Gamma}^{\alpha}_{\mu\nu}$:
\begin{eqnarray*}
\mathbf{\Gamma}^{\alpha}_{\mu\nu}-\mathbf{\Gamma}^{\alpha}_{\nu\mu}&=&
({\Gamma}^{\alpha}_{\mu\nu}-\frac{1}{(n-1)}\delta^{\alpha}_{\mu}C_{\nu})
-({\Gamma}^{\alpha}_{\nu\mu}-\frac{1}{(n-1)}\delta^{\alpha}_{\nu}C_{\mu})\\
&=&{\Lambda}^{\alpha}_{\mu\nu}-\frac{1}{(n-1)}\set{\delta^{\alpha}_{\mu}C_{\nu}-\delta^{\alpha}_{\nu}C_{\mu}}=
T^{\alpha}_{\mu\nu}.
\end{eqnarray*}

\begin{eqnarray*}
\mathfrak{U}_{\nu\rho}\{\mathbf{\Gamma}^{\alpha}_{\mu\sigma,\nu}+
\mathbf{\Gamma}^{\epsilon}_{\mu\sigma}\mathbf{\Gamma}^{\alpha}_{\epsilon\nu}\}
 &=& ({\Gamma}^{\alpha}_{\mu\sigma}-\frac{1}{(n-1)}\,\delta^{\alpha}_{\mu}C_{\sigma})_{\!,\nu}-
({\Gamma}^{\alpha}_{\mu\nu}-\frac{1}{(n-1)}\,\delta^{\alpha}_{\mu}C_{\nu})_{\!,\sigma}\\
&&+({\Gamma}^{\epsilon}_{\mu\sigma}-\frac{1}{(n-1)}\,\delta^{\epsilon}_{\mu}C_{\sigma})
({\Gamma}^{\alpha}_{\epsilon\nu}-\frac{1}{(n-1)}\,\delta^{\alpha}_{\epsilon}C_{\nu})\\
&&-({\Gamma}^{\epsilon}_{\mu\nu}-\frac{1}{(n-1)}\,\delta^{\epsilon}_{\mu}C_{\nu})
({\Gamma}^{\alpha}_{\epsilon\sigma}-\frac{1}{(n-1)}\,\delta^{\alpha}_{\epsilon}C_{\sigma})\\
&=& {\Gamma}^{\alpha}_{\mu\sigma,\nu}-{\Gamma}^{\alpha}_{\mu\nu,\sigma}+
{\Gamma}^{\epsilon}_{\mu\sigma}{\Gamma}^{\alpha}_{\epsilon\nu}-
{\Gamma}^{\epsilon}_{\mu\nu}{\Gamma}^{\alpha}_{\epsilon\sigma}\\
&& -\frac{1}{(n-1)}\delta^{\alpha}_{\mu}C_{\sigma,\nu}+\frac{1}{(n-1)}\,\delta^{\alpha}_{\mu}C_{\nu,\sigma}\\
&&-\frac{1}{(n-1)}\Gamma^{\alpha}_{\mu\sigma}C_{\nu}-\frac{1}{(n-1)}\Gamma^{\alpha}_{\mu\nu}C_{\sigma}
+\frac{1}{(n-1)^{2}}\,\delta^{\alpha}_{\mu}C_{\sigma}C_{\nu}\\
&&+\frac{1}{(n-1)}\Gamma^{\alpha}_{\mu\nu}C_{\sigma}+\frac{1}{(n-1)}\Gamma^{\alpha}_{\mu\sigma}C_{\nu}
-\frac{1}{(n-1)^{2}}\,\delta^{\alpha}_{\mu}C_{\nu}C_{\sigma}\\
&=&{R}^{\alpha}_{\mu\nu\sigma}-\frac{1}{(n-1)}\set{\delta^{\alpha}_{\mu}C_{\sigma,\nu}-\delta^{\alpha}_{\mu}C_{\nu,\sigma}}\\
&=&\frac{1}{(n-1)}\set{\delta^{\alpha}_{\mu}C_{\nu,\sigma}-\delta^{\alpha}_{\mu}C_{\sigma,\nu}}=
K^{\alpha}_{\mu\nu\sigma},
\end{eqnarray*}
since the curvature tensor of the Weitzenb\"ock connection vanishes identically.\\
This completes the proof.
 \end{proof}

\vspace{5pt}
Let $\widehat{\Gamma}^\alpha_{\mu\nu}$ be the symmetric part of the Weitzenb\"ock connection: $\widehat{\Gamma}^\alpha_{\mu\nu}=\frac{1}{2}(\Gamma^\alpha_{\mu\nu}+\Gamma^\alpha_{\nu\mu})$. This is a symmetric connection with curvature $\widehat{R}^\alpha_{\mu\nu\sigma}$ \cite{Sid-Ahmed}.
Under the conformal change (\ref{eq.7}), one can show that:
\begin{equation}\label{SymCon}
    \overline{\widehat{\Gamma}}\,^{\alpha}_{\mu\nu}=\widehat{\Gamma}^{\alpha}_{\mu\nu}+
  \frac{1}{2}(\delta^{\alpha}_{\mu}\rho_{\nu}+\delta^{\alpha}_{\nu}\rho_{\mu}),
\end{equation}

\begin{equation}\label{curvSymCon}
    \overline{\widehat{R}}\,^{\!\alpha}_{\mu\nu\sigma}=\widehat{R}^{\alpha}_{\mu\nu\sigma}+ \frac{1}{2}\mathfrak{U}_{\nu\sigma}
  \set{\delta^{\alpha}_{\sigma}\rho_{\mu\widehat{\mid}\nu}
  +\frac{1}{2}\delta^{\alpha}_{\nu}\rho_{\sigma}\rho_{\mu}},
\end{equation}
where $\widehat{|}$ denotes the covariant derivatives with respect to $\widehat{\Gamma}^\alpha_{\mu\nu}$.

\vspace{8pt}
 \noindent\textbf{Theorem B.}
\emph{Let $(M,\stackrel[i]{}{\lambda})$ be  an AP-space of dimension $n\geq2$. The tensor
\begin{eqnarray*}\label{e.10ss}
B^{\alpha}_{\mu\nu\sigma}&:=&\frac{1}{4}\,\mathfrak{U}_{\nu,\sigma}
\set{2\Lambda^{\alpha}_{\mu\nu|\sigma}+\Lambda^{\epsilon}_{\mu\nu}\Lambda^{\alpha}_{\sigma\epsilon}
+\Lambda^{\epsilon}_{\sigma\nu}\Lambda^{\alpha}_{\epsilon\mu}}\\
&&-\frac{1}{2(n-1)}\mathfrak{U}_{\nu,\sigma}\set{\delta^{\alpha}_{\mu}C_{\sigma,\nu}+\delta^{\alpha}_{\sigma}C_{\mu\widehat{|}\nu}
-\frac{1}{2(n-1)}\delta^{\alpha}_{\nu}C_{\mu}C_{\sigma}}
\end{eqnarray*}
is conformally invariant. Moreover,  $B^{\alpha}_{\mu\nu\sigma}$
is precisely the curvature tensor of a conformal connection on $M$.}

 \begin{proof} Let the covariant derivative with respect to $\overline{\widehat{\Gamma}}\,^{\alpha}_{\mu\nu}$ be denoted by $\widehat{||}$. Using Equations (\ref{C}) and (\ref{SymCon}), one can show that
 \begin{eqnarray*}\label{C1}
    \overline{C}_{\sigma,\nu}&=&C_{\sigma,\nu}+(n-1)\rho_{\sigma,\nu}\\
  \overline{C}_{\mu\widehat{\|}\nu}&=&C_{\mu\widehat{|}\nu}+(n-1)\rho_{\mu\widehat{|}\nu}-
  \frac{1}{2}(C_{\mu}\rho_{\nu}+C_{\nu}\rho_{\mu})-(n-1)\rho_{\mu}\rho_{\nu},
 \end{eqnarray*}

 We show that the tensor $B^{\alpha}_{\mu\nu\sigma}$ is conformally invariant.
   From the above two relations together with  (\ref{C}), (\ref{curvSymCon}) and Theorem 1\textbf{(b)} of \cite{Sid-Ahmed}, we get, after some manipulations,

 \begin{eqnarray*}
\overline{B}^{\,\alpha}_{\mu\nu\sigma}&=&{\overline{\widehat{R}}}\,^{\,\alpha}_{\mu\nu\sigma}-\frac{1}{2(n-1)}
\{\delta^{\alpha}_{\mu}\overline{C}_{\sigma,\nu}-\delta^{\alpha}_{\mu}\overline{C}_{\nu,\sigma}
+\delta^{\alpha}_{\sigma}\overline{C}_{\mu\widehat{\|}\nu}-\delta^{\alpha}_{\nu}\overline{C}_{\mu\widehat{\|}\sigma}\\
&&-\frac{1}{2(n-1)}\,\delta^{\alpha}_{\nu}\overline{C}_{\mu}\overline{C}_{\sigma}
+\frac{1}{2(n-1)}\,\delta^{\alpha}_{\sigma}\overline{C}_{\mu}\overline{C}_{\nu}\}\\
&=&\widehat{R}^{\alpha}_{\mu\nu\sigma}+ \frac{1}{2}
 \{ \delta^{\alpha}_{\sigma}\rho_{\mu\widehat{\mid}\nu}-\delta^{\alpha}_{\nu}\rho_{\mu\widehat{\mid}\sigma}
  +\frac{1}{2}\,\delta^{\alpha}_{\nu}\rho_{\sigma}\rho_{\mu}-\frac{1}{2}\,\delta^{\alpha}_{\sigma}\rho_{\nu}\rho_{\mu}\}\\
 && -\frac{1}{2(n-1)}
\{\delta^{\alpha}_{\mu}(C_{\sigma,\nu}+(n-1)\rho_{\sigma,\nu})-\delta^{\alpha}_{\mu}(C_{\nu,\sigma}+(n-1)\rho_{\nu,\sigma})\\
&&+\delta^{\alpha}_{\sigma}(C_{\mu\widehat{|}\nu}+(n-1)\rho_{\mu\widehat{|}\nu}-
  \frac{1}{2}(C_{\mu}\rho_{\nu}+C_{\nu}\rho_{\mu})-(n-1)\rho_{\mu}\rho_{\nu})\\
  &&-\delta^{\alpha}_{\nu}(C_{\mu\widehat{|}\sigma}+(n-1)\rho_{\mu\widehat{|}\sigma}-
  \frac{1}{2}(C_{\mu}\rho_{\sigma}+C_{\sigma}\rho_{\mu})-(n-1)\rho_{\mu}\rho_{\sigma})\\
&&-\frac{1}{2(n-1)}\,\delta^{\alpha}_{\nu}(C_{\mu}+(n-1)\rho_{\mu})(C_{\sigma}+(n-1)\rho_{\sigma})\\
&&+\frac{1}{2(n-1)}\,\delta^{\alpha}_{\sigma}(C_{\mu}+(n-1)\rho_{\mu})(C_{\nu}+(n-1)\rho_{\nu})\}\\
&=&{{\widehat{R}}}^{\,\alpha}_{\mu\nu\sigma}-\frac{1}{2(n-1)}
\{\delta^{\alpha}_{\mu}{C}_{\sigma,\nu}-\delta^{\alpha}_{\mu}{C}_{\nu,\sigma}
+\delta^{\alpha}_{\sigma}{C}_{\mu\widehat{|}\nu}-\delta^{\alpha}_{\nu}{C}_{\mu\widehat{|}\sigma}\\
&&-\frac{1}{2(n-1)}\,\delta^{\alpha}_{\nu}{C}_{\mu}{C}_{\sigma}
+\frac{1}{2(n-1)}\,\delta^{\alpha}_{\sigma}{C}_{\mu}{C}_{\nu}\}=
B^{\alpha}_{\mu\nu\sigma}.
 \end{eqnarray*}

Now, define the connection:
\begin{equation}\label{e.9sy}
\mathbf{\widehat{\Gamma}}^{\alpha}_{\mu\nu}:={\widehat{\Gamma}}^{\alpha}_{\mu\nu}-\frac{1}{2(n-1)}
(\delta^{\alpha}_{\mu}C_{\nu}+\delta^{\alpha}_{\nu}C_{\mu}).
\end{equation}
One can easily show that this connection is conformal. On the other hand, we prove that the tensor  $B^{\alpha}_{\mu\nu\sigma}$
is the curvature tensor of the conformal connection $\mathbf{\widehat{\Gamma}}^{\alpha}_{\mu\nu}$:

\begin{eqnarray*}
&~&\mathfrak{U}_{\nu,\sigma}\set{\mathbf{\widehat{\Gamma}}^{\alpha}_{\mu\sigma,\nu}+
\mathbf{\widehat{\Gamma}}^{\epsilon}_{\mu\sigma}\mathbf{\widehat{\Gamma}}^{\alpha}_{\epsilon\nu}}\\
&=& ({\widehat{\Gamma}}^{\alpha}_{\mu\sigma}-\frac{1}{2(n-1)}(\delta^{\alpha}_{\mu}C_{\sigma}+\delta^{\alpha}_{\sigma}C_{\mu}))_{\!,\nu}-
({\widehat{\Gamma}}^{\alpha}_{\mu\nu}-\frac{1}{2(n-1)}(\delta^{\alpha}_{\mu}C_{\nu}+\delta^{\alpha}_{\nu}C_{\mu}))_{\!,\sigma}\\
&&+({\widehat{\Gamma}}^{\epsilon}_{\mu\sigma}-\frac{1}{2(n-1)}(\delta^{\epsilon}_{\mu}C_{\sigma}+\delta^{\epsilon}_{\sigma}C_{\mu}))
({\widehat{\Gamma}}^{\alpha}_{\epsilon\nu}-\frac{1}{2(n-1)}(\delta^{\alpha}_{\epsilon}C_{\nu}+\delta^{\alpha}_{\nu}C_{\epsilon}))\\
&&-({\widehat{\Gamma}}^{\epsilon}_{\mu\nu}-\frac{1}{2(n-1)}(\delta^{\epsilon}_{\mu}C_{\nu}+\delta^{\epsilon}_{\nu}C_{\mu}))
({\widehat{\Gamma}}^{\alpha}_{\epsilon\sigma}-\frac{1}{2(n-1)}(\delta^{\alpha}_{\epsilon}C_{\sigma}+\delta^{\alpha}_{\sigma}C_{\epsilon}))\\
&=& {\widehat{\Gamma}}^{\alpha}_{\mu\sigma,\nu}-{\widehat{\Gamma}}^{\alpha}_{\mu\nu,\sigma}+
{\widehat{\Gamma}}^{\epsilon}_{\mu\sigma}{\widehat{\Gamma}}^{\alpha}_{\epsilon\nu}-
{\widehat{\Gamma}}^{\epsilon}_{\mu\nu}{\widehat{\Gamma}}^{\alpha}_{\epsilon\sigma}\\
&& -\frac{1}{2(n-1)}(\delta^{\alpha}_{\mu}C_{\sigma,\nu}+\delta^{\alpha}_{\sigma}C_{\mu,\nu})
+\frac{1}{2(n-1)}(\delta^{\alpha}_{\mu}C_{\nu,\sigma}+\delta^{\alpha}_{\nu}C_{\mu,\sigma})\\
&&-\frac{1}{2(n-1)}(\widehat{\Gamma}^{\alpha}_{\mu\sigma}C_{\nu}+\delta^{\alpha}_{\nu}\widehat{\Gamma}^{\epsilon}_{\mu\sigma}C_{\epsilon})
-\frac{1}{2(n-1)}(\widehat{\Gamma}^{\alpha}_{\mu\nu}C_{\sigma}+\widehat{\Gamma}^{\alpha}_{\sigma\nu}C_{\mu})\\
&&+\frac{1}{4(n-1)^{2}}(\delta^{\alpha}_{\mu}C_{\nu}C_{\sigma}+\delta^{\alpha}_{\nu}C_{\mu}C_{\sigma}+
\delta^{\alpha}_{\sigma}C_{\nu}C_{\mu}+\delta^{\alpha}_{\nu}C_{\sigma}C_{\mu})\\
&&+\frac{1}{2(n-1)}(\widehat{\Gamma}^{\alpha}_{\mu\nu}C_{\sigma}+\delta^{\alpha}_{\sigma}\widehat{\Gamma}^{\epsilon}_{\mu\nu}C_{\epsilon})
+\frac{1}{2(n-1)}(\widehat{\Gamma}^{\alpha}_{\mu\sigma}C_{\nu}+\widehat{\Gamma}^{\alpha}_{\nu\sigma}C_{\mu})\\
&&-\frac{1}{4(n-1)^{2}}(\delta^{\alpha}_{\mu}C_{\sigma}C_{\nu}+\delta^{\alpha}_{\sigma}C_{\mu}C_{\nu}+
\delta^{\alpha}_{\nu}C_{\sigma}C_{\mu}+\delta^{\alpha}_{\sigma}C_{\nu}C_{\mu})\\
&=&{{\widehat{R}}}^{\,\alpha}_{\mu\nu\sigma}-\frac{1}{2(n-1)}
\{\delta^{\alpha}_{\mu}{C}_{\sigma,\nu}-\delta^{\alpha}_{\mu}{C}_{\nu,\sigma}
+\delta^{\alpha}_{\sigma}{C}_{\mu\widehat{|}\nu}-\delta^{\alpha}_{\nu}{C}_{\mu\widehat{|}\sigma}\\
&&-\frac{1}{2(n-1)}\,\delta^{\alpha}_{\nu}{C}_{\mu}{C}_{\sigma}
+\frac{1}{2(n-1)}\,\delta^{\alpha}_{\sigma}{C}_{\mu}{C}_{\nu}\}=
B^{\alpha}_{\mu\nu\sigma}.
\end{eqnarray*}
This completes the proof.
 \end{proof}

The next lemma is used for proving Theorem C below.
\begin{lem}\label{L.1}
Under the conformal change $(\ref{eq.7})$, we have
\begin{description}
  \item[(a)]$\overline{C}_{\sigma}=C_{\sigma}+(n-1)\rho_{\sigma}\,$,  $\overline{C}^{\sigma}=e^{-2\rho(x)}(C^{\sigma}+(n-1)\rho^{\sigma})$

  \item[(b)]$\overline{C}^{2}=e^{-2\rho(x)}(C^{2}+2(n-1)C_{\epsilon}\rho^{\epsilon}+(n-1)^{2}\rho^{2})$

  \item[(c)]$\overline{C}_{\mu;;\nu}=C_{\mu;\nu}+(n-1)\rho_{\mu;\nu}
  -(C_{\mu}\rho_{\nu}+C_{\nu}\rho_{\mu}-g_{\mu\nu}C_{\epsilon}\rho^{\epsilon})-(n-1)(2\rho_{\mu}\rho_{\nu}-g_{\mu\nu}\rho^{2})$

  \item[(d)]$\overline{C}^{\alpha}_{\,\,;;\nu}=e^{-2\rho(x)}\{C^{\alpha}_{\,\,;\nu}+(n-1)\rho^{\alpha}_{\,\,;\nu}
  +(\delta^{\alpha}_{\nu}C^{\epsilon}\rho_{\epsilon}-C^{\alpha}\rho_{\nu}-C_{\nu}\rho^{\alpha})
  +(n-1)(\delta^{\alpha}_{\nu}\rho^{2}-2\rho^{\alpha} \rho_{\nu})\}$,
 \end{description}
where $C^{\alpha}:=g^{\alpha\epsilon}C_{\epsilon}, \,\, \overline{C}^{2}:=\overline{C}_{\epsilon}\overline{C}^{\epsilon}$ and $;$ and $;;$ are the covariant derivatives with respect to $\overcirc{\Gamma}\,^{\alpha}_{\mu\nu}$ and\,    $\overline{\overcirc{\Gamma}}\,^{\alpha}_{\mu\nu}$, respectively.
\end{lem}

 \vspace{3pt}

\noindent\textbf{Theorem C.}
\emph{Let $(M,\stackrel[i]{}{\lambda})$ be  an AP-space of dimension $n\geq2$. The tensor
\begin{eqnarray}\label{e.10rr}
  Q^{\alpha}_{\mu\nu\sigma}&:=&\mathfrak{U}_{\nu,\sigma}\{\gamma^{\alpha}_{\mu\nu|\sigma} +
  \gamma^{\epsilon}_{\mu\sigma}\gamma^{\alpha}_{\epsilon\nu}+\frac{1}{2}\gamma^{\alpha}_{\mu\epsilon}\Lambda^{\epsilon}_{\nu\sigma}\}-\frac{1}{(n-1)}
\mathfrak{U}_{\nu,\sigma}\{\delta^{\alpha}_{\mu}C_{\sigma,\nu}+\delta^{\alpha}_{\sigma}C_{\mu;\nu}+
g_{\mu\sigma}C^{\alpha}_{\,\,\,;\nu} \nonumber\\
&&- \frac{1}{(n-1)}(\delta^{\alpha}_{\nu}C_{\mu}C_{\sigma} -\delta^{\alpha}_{\nu}g_{\mu\sigma}C^{2}+g_{\mu\sigma}C_{\nu}C^{\alpha})   \}
\end{eqnarray}
is conformally invariant. Moreover,  $Q^{\alpha}_{\mu\nu\sigma}$
is precisely the  curvature tensor of a conformal connection on $M$. }

\vspace{5pt}
  In fact, the proof is similar to that of the preceding theorem, taking into account Lemma \ref{L.1} above and Theorem 1\textbf{(c)} of \cite{Sid-Ahmed}.
 The conformal connection whose curvature tensor coincides with $Q^{\alpha}_{\mu\nu\sigma}$ is given by
 $${\bf\stackrel[]{o}\Gamma}\,^{\alpha}_{\mu\nu}:=\,\overcirc{\Gamma}\,^{\alpha}_{\mu\nu}
-\frac{1}{(n-1)}(\delta^{\alpha}_{\mu}C_{\nu}+\delta^{\alpha}_{\nu}C_{\mu}-g_{\mu\nu}C^{\alpha}).$$


\providecommand{\bysame}{\leavevmode\hbox
to3em{\hrulefill}\thinspace}
\providecommand{\MR}{\relax\ifhmode\unskip\space\fi MR }
\providecommand{\MRhref}[2]{%
  \href{http://www.ams.org/mathscinet-getitem?mr=#1}{#2}
} \providecommand{\href}[2]{#2}


\begin{thebibliography}{10}

\bibitem{Abed}
S.~H. Abed, \emph{Conformal $\beta$-changes in \textsc{F}insler
spaces}, Proc. Math. Phys. Soc. Egypt, \textbf{86} (2008), 79-89.
arXiv: math. DG/0602404.

\bibitem{Br.1}
F. Brickell and R. S. Clark, Differentiable manifolds, Van Nostrand Reinhold
Co., 1970.

\bibitem{WN} W. El Hanafy and G. G. L. Nashed, \emph{Exact teleparallel gravity of binary black holes}, Astrophys. Space Sci., \textbf{361} (2016), 68. DOI: 10.1007/s10509-016-2662-y. arXiv: 1507.07377 [gr-qc].

\bibitem{Hashiguchi}
M.~Hashiguchi, \emph{On conformal transformation of {F}insler metrics},
  J. Math. Kyoto Univ., {\textbf{16}} (1976), 25-50.

\bibitem{Izumi1}
H.~Izumi, \emph{Conformal transformations of {F}insler spaces}, I,
Tensor, N. S., {\textbf{31}} (1977), 33-41.

\bibitem{Izumi2}
H.~Izumi, \emph{Conformal transformations of \textsc{F}insler
spaces}, II, Tensor, N. S., \textbf{34} (1980), 337-359.

\bibitem{Knebelman}
M. S. Knebelman: {\it Conformal geometry of generalized metric
spaces}, Proc. Nat. Acad. Sci. USA, \textbf{15} (1929), 376-379.

\bibitem{Matsumoto}
 M. Matsumoto, \emph{
Conformally closed Finsler spaces}, Balkan J. Geom. Appl.,
\textbf{4}, 1 (1999), 117-128.

\bibitem{NW} G. G. L. Nashed and W. El Hanafy, \emph{A Built-in Inflation in the $f(T)$-Cosmology}, Europ. Phys. J. C, \textbf{74 }(2014), 3099. DOI: 10.1140/epjc/s10052-014-3099-5. arXiv: 1403.0913 [gr-qc].

\bibitem{Obata}
M. Obata,  {\it Conformal transformation in Riemannian manifolds}, Sugaku Math. Soc. Japan, \textbf{14} (1963), 152-164.

\bibitem{Rob.1}
 H. P. Robertson, \emph{Groups of motion in spaces admitting absolute parallelism},
Ann. Math, Princeton (2), \textbf{33} (1932), 496-520.

\bibitem{Shirafuji} T. Shirafuji, G. G. Nashed and Y. Kobayashi \emph{Equivalence principle in the new general relativity}, Prog. Theoret. Phys., \textbf{96} (1996), 933-947. DOI: 10.1143/PTP.96.933. arXiv: gr-qc/9609060.

\bibitem{Wanas} M. I. Wanas, \emph{Geometrical structures for cosmological applications}, Astrophys. Space Sci., \textbf{127} (1986), 21-25. DOI: 10.1007/BF00637758.

\bibitem{Wan.2}
M. I. Wanas, \emph{Absolute parallelism geometry}: Developments, applications and
problems, Stud. Cercet, Stiin. Ser. Mat. Univ. Bacau, No. 10 (2001), 297-309. arXiv:gr-qc/0209050.

\bibitem{charge} M. I. Wanas,\emph{ On the relation between mass and charge: A pure geometric approach},
    Int. J. Geom. Meth. Mod. Phys., \textbf{4} (2007), 373-388. arXiv: gr-qc/0703036.

\bibitem{WYS} M. I. Wanas, N. L. Youssef and A. M. Sid-Ahmed, \emph{Teleparallel Lagrange geometry and a unified field theory}, Class. Quantum Grav., \textbf{27} (2010) 045005 (29pp). Doi: 10.1088/0264-9381/27/4/045005. arXiv: 0905.0209  [gr-qc].

\bibitem{Yan.1}
K. Yano, {Integral formulas in Riemannian geometry}, Marcel Dekker, INC., New York, 1970.

\bibitem{NAS1}
Nabil~L. Youssef, S.~H. Abed and A.~Soleiman, \emph{A global theory
of conformal {F}insler geometry}, Tensor, N. S., {\textbf{69}}
(2008), 155--178. arXiv: math. DG/0610052.

\bibitem{NAS2}
Nabil L. Youssef, S. H. Abed and A. Soleiman,
\emph{Conformal change of special Finsler spaces}, Balkan J. Geom. Appl., \textbf{15} (2010), 146-158. arXiv: 0908.0696 [math.DG].

\bibitem{Waleed}
Nabil L. Youssef and Waleed A. Elsayed, \emph{A global approach to absolute parallelism geometry}, Rep. Math. Phys., \textbf{72} (2013), 1-23. arXiv: 1209.1379 [gr-qc].

\bibitem{Sid-Ahmed}
Nabil L. Youssef and Amr M. Sid-Ahmed, \emph{Linear connections and curvature tensors
 in the geometry of parallizable manifold}, Rep. Math. Phys., \textbf{60} (2007), 39-53. arXiv: gr-qc/0604111.




\end{thebibliography}
\end{document}